\documentclass[10pt,reqno]{amsart}
\usepackage{enumerate, latexsym, amsmath, amsfonts, amssymb, amsthm, color}
\def\pmod #1{\ ({\rm{mod}}\ #1)}
\def\Z{\Bbb Z}
\def\N{\Bbb N}

\def\Q{\Bbb Q}

\def\l{\left}
\def\r{\right}
\def\bg{\bigg}
\def\({\bg(}
\def\){\bg)}
\def\t{\text}
\def\f{\frac}

\def\mo{{\rm{mod}\ }}

\def\ls{\leqslant}
\def\gs{\geqslant}

\def\bi{\binom}
\def\al{\alpha}

\def\eq{\equiv}

\def\da{\delta}

\def\Proof{\noindent{\it Proof}}
\def\Ack{\medskip\noindent {\bf Acknowledgments}}
\theoremstyle{plain}
\newtheorem{theorem}{Theorem}

\newtheorem{lemma}[theorem]{Lemma}

\newtheorem{conjecture}{Conjecture}
\theoremstyle{definition}

\theoremstyle{remark}
\newtheorem{remark}[theorem]{Remark}

 \vspace{4mm}

\begin{document}

\hbox{Int. J. Number Theory 16 (2020), no.\,8, 1833--1858.}
\medskip

\title
[{Quadratic residues and quartic residues modulo primes}]
{Quadratic residues and quartic residues \\ modulo primes}

\author[Z.-W. Sun] {Zhi-Wei Sun}

\address{Department of Mathematics, Nanjing
University, Nanjing 210093, People's Republic of China}
\email{zwsun@nju.edu.cn}

\keywords{Quadratic residues, quartic residues, congruences, Lucas sequences.
\newline \indent 2020 {\it Mathematics Subject Classification}.
Primary 11A15; Secondary 11A07, 11B39.
}

\begin{abstract}
In this paper we study some products related to quadratic residues and quartic residues modulo primes.
Let $p$ be an odd prime and let $A$ be any integer. We determine completely the product
$$f_p(A):=\prod_{1\ls i,j\ls(p-1)/2\atop p\nmid i^2-Aij-j^2}(i^2-Aij-j^2)$$ modulo $p$; for example, if $p\eq1\pmod4$ then
$$f_p(A)
\equiv\begin{cases}-(A^2+4)^{(p-1)/4}\pmod p&\text{if}\ (\frac{A^2+4}p)=1,
\\(-A^2-4)^{(p-1)/4}\pmod p&\text{if}\ (\frac{A^2+4}p)=-1,\end{cases}$$
where $(\frac{\cdot}p)$ denotes the Legendre symbol. We also determine
$$\prod^{(p-1)/2}_{i,j=1\atop p\nmid 2i^2+5ij+2j^2}\l(2i^2+5ij+2j^2\r)
\ \text{and}\ \prod^{(p-1)/2}_{i,j=1\atop p\nmid 2i^2-5ij+2j^2}\l(2i^2-5ij+2j^2\r)$$
modulo $p$.
\end{abstract}
\maketitle

\section{Introduction}
\setcounter{theorem}{0}
\setcounter{corollary}{0}
\setcounter{equation}{0}

For $n\in\Z^+=\{1,2,3,\ldots\}$ and $x=a/b$ with $a,b\in\Z$, $b\not=0$ and $\gcd(b,n)=1$,
we let $\{x\}_n$ denote the unique integer $r\in\{0,\ldots,n-1\}$ with $r\eq x\pmod n$ (i.e., $a\eq br\pmod n$). The well-known Gauss Lemma (see, e.g., \cite[p.\,52]{IR}) states that for any odd prime $p$ and integer $x\not\eq0\pmod p$
we have
\begin{equation}\label{1.1}\l(\f xp\r)=(-1)^{|\{1\ls k<p/2:\ \{kx\}_p>p/2\}|},
\end{equation}
where $(\f{\cdot}p)$ is the Legendre symbol. This was extended to Jacobi symbols by M. Jenkins \cite{J} in 1867, who showed (by an elementary method) that for any positive odd integer $n$ and integer $x$ with $\gcd(x,n)=1$ we have
\begin{equation}\label{1.2}\l(\f xn\r)=(-1)^{|\{1\ls k<n/2:\ \{kx\}_n>n/2\}|},
\end{equation}
where $(\f{\cdot}n)$ is the Jacobi symbol.  In the textbook \cite[Chapters 11-12]{R}, H. Rademacher supplied a proof of Jenkins' result by using subtle properties of quadratic Gauss sums.

Now we present our first new theorem.

\begin{theorem}\label{Th1.1} Let $n$ be a positive odd integer, and let $x\in\Z$ with $\gcd(x(1-x),n)=1$. Then
\begin{equation}\label{1.3}(-1)^{|\{1\ls k<n/2:\ \{kx\}_n>k\}|}=\l(\f{2x(1-x)}n\r).
\end{equation}
Also,
\begin{align}\label{1.4}(-1)^{|\{1\ls k<n/2:\ \{kx\}_n>n/2\ \&\ \{k(1-x)\}_n>n/2\}|}=&\l(\f 2n\r),
\\\label{1.5}(-1)^{|\{1\ls k<n/2:\ \{kx\}_n<n/2\ \&\ \{k(1-x)\}_n<n/2\}|}=&\l(\f {2x(x-1)}n\r),
\end{align}
and
\begin{equation}\label{1.6}(-1)^{|\{1\ls k<n/2:\ \{kx\}_n>n/2>\{k(1-x)\}_n\}|}=\l(\f {2x}n\r).
\end{equation}
\end{theorem}

 Let $p$ be an odd prime, and let $a,b,c\in\Z$ and
\begin{equation}\label{1.7}S_p(a,b,c):=\prod_{1\ls i<j\ls p-1\atop p\nmid ai^2+bij+cj^2}(ai^2+bij+cj^2).
\end{equation}
Using Theorem \ref{Th1.1}, together with \cite[Theorem 1.2]{S19b}, we completely determine $S_p(a,b,c)$ mod $p$ in terms of Legendre symbols.

\begin{theorem}\label{Th1.2} Let $p$ be an odd prime, and let $a,b,c\in\Z$ and $\Delta=b^2-4ac$.
When $p\nmid ac(a+b+c)$, we have
\begin{equation}\label{1.8}S_p(a,b,c)\eq
\begin{cases}(\f{a(a+b+c)}p)\pmod p&\t{if}\ p\mid \Delta,\\-(\f{ac(a+b+c)\Delta}p)\pmod p&\t{if}\ p\nmid \Delta.\end{cases}
\end{equation}
In the case $p\mid ac(a+b+c)$, we have
\begin{equation}\label{1.9}
S_p(a,b,c)\eq\begin{cases}0\pmod p&\t{if}\ p\mid a,\ p\mid b\ \t{and}\ p\mid c,
\\-(\f {-a}p)\pmod p&\t{if}\ p\nmid a,\ p\mid b\ \t{and}\ p\mid c,
\\-(\f bp)\pmod p&\t{if}\ p\mid a,\ p\nmid b\ \t{and}\ p\mid c,
\\-(\f {-c}p)\pmod p&\t{if}\ p\mid a,\ p\mid b\ \t{and}\ p\nmid c,
\\-(\f cp)\pmod p&\t{if}\ p\mid a,\ p\nmid bc\ \t{and}\ p\mid b+c,
\\-(\f {a}p)\pmod p&\t{if}\ p\nmid ab,\ p\mid a+b\ \t{and}\ p\mid c,
\\-(\f {-a}p)\pmod p&\t{if}\ p\nmid ac,\ p\mid a-c \ \t{and}\ p\mid a+b+c,
\\(\f{-ac}p)\pmod p&\t{if}\ p\nmid ac(a-c)\ \t{and}\ p\mid a+b+c,
\\(\f{-a(a+b)}p)\pmod p&\t{if}\ p\nmid ab(a+b)\ \t{and}\ p\mid c,
\\(\f{-c(b+c)}p)\pmod p&\t{if}\ p\mid a\ \t{and}\ p\nmid bc(b+c),
\end{cases}\end{equation}
\end{theorem}

We will prove Theorem \ref{Th1.1} and those parts of Theorem \ref{Th1.2} not covered by \cite[Theorem 1.2]{S19b} in Section 2.

Let $p$ be an odd prime. For $a,b,c\in\Z$ we introduce
\begin{equation}\label{1.10} T_p(a,b,c):=\prod_{i,j=1\atop p\nmid ai^2+bij+cj^2}^{(p-1)/2}(ai^2+bij+cj^2).\end{equation}
Our following theorem determines $T_p(1,-(a+b),-1)$ modulo $p$ for all $a,b\in\Z$ with $ab\eq-1\pmod p$.

\begin{theorem}\label{Th1.3} Let $p$ be any odd prime, and let $a,b\in\Z$ with $ab\eq-1\pmod p$.
Set
\begin{equation}\label{1.11} \{a,b\}_p:=\prod^{(p-1)/2}_{i,j=1\atop i\not\eq aj,bj\pmod p}(i-aj)(i-bj),
\end{equation}
which is congruent to $T_p(1,-(a+b),-1)$ modulo $p$.

{\rm (i)} We have
\begin{equation}\label{1.12}-\{a,b\}_p\eq\begin{cases}(\f{a-b}p)\pmod p&\t{if}\ p\eq1\pmod4\ \&\ p\nmid (a-b),
\\(\f{a(a-b)}p)=(\f{a^2+1}p)\pmod p&\t{if}\ p\eq3\pmod4.\end{cases}
\end{equation}

{\rm (ii)} If $a\eq b\pmod p$ and $p\eq1\pmod 8$, then
$$\{a,b\}_p\eq(-1)^{(p+7)/8}\f{p-1}2!\pmod p.$$
If $p\eq5\pmod 8$ and $a\eq b\eq(-1)^k((p-1)/2)!\pmod p$ with $k\in\{0,1\}$, then
$$\{a,b\}_p\eq (-1)^{k+(p-5)/8}\pmod p.$$
\end{theorem}

Our proof of Theorem \ref{Th1.3} will be given in Section 3.

For any $A\in\Z$, we define the Lucas sequences $\{u_n(A)\}_{n\gs0}$ and $\{v_n(A)\}_{n\gs0}$ by
$$u_0(A)=0,\ u_1(A)=1,\ \t{and}\ u_{n+1}(A)=Au_n(A)+u_{n-1}(A)\ \t{for}\ n=1,2,3,\ldots,$$
and
$$v_0(A)=2,\ v_1(A)=A,\ \t{and}\ v_{n+1}(A)=Av_n(A)+v_{n-1}(A)\ \t{for}\ n=1,2,3,\ldots.$$
It is well known that
$$u_n(A)=\f{\al^n-\beta^n}{\al-\beta}\ \ \t{and}\ \ v_n(A)=\al^n+\beta^n$$
for all $n\in\N=\{0,1,2,\ldots\},$ where
$$\al=\f{A+\sqrt{A^2+4}}2\ \ \t{and}\ \ \beta=\f{A-\sqrt{A^2+4}}2.$$
Thus
\begin{equation}\label{1.13}
\l(\f{A\pm\sqrt{A^2+4}}2\r)^n=\f{v_n(A)\pm u_n(A)\sqrt{A^2+4}}2\quad\t{for all}\ n\in\N.
\end{equation}

Now we state our fourth theorem which determines $T_p(1,-A,-1)$ for any odd prime $p$ and integer $A$.

\begin{theorem}\label{Th1.4} Let $p$ be an odd prime and let $A\in\Z$.

{\rm (i)} Suppose that
$p\mid(A^2+4)$. Then $p\eq1\pmod4$, $A/2\eq (-1)^k((p-1)/2)!\pmod p$ for some $k\in\{0,1\}$, and
\begin{equation}\label{1.14}T_p(1,-A,-1)\eq
\begin{cases}(-1)^{(p+7)/8}((p-1)/2)!\pmod p&\t{if}\ p\eq1\pmod 8,
\\(-1)^{k+(p-5)/8}\pmod p&\t{if}\ p\eq5\pmod 8.\end{cases}
\end{equation}

{\rm (ii)} When $(\f{A^2+4}p)=1$, we have
\begin{equation}\label{1.15}T_p(1,-A,-1)
\eq\begin{cases}-(A^2+4)^{(p-1)/4}\pmod p&\t{if}\ p\eq1\pmod4,
\\-(A^2+4)^{(p+1)/4}u_{(p-1)/2}(A)/2\pmod p&\t{if}\ p\eq3\pmod4.\end{cases}
\end{equation}

{\rm (iii)} When $(\f{A^2+4}p)=-1$, we have
\begin{equation}\label{1.16}T_p(1,-A,-1)
\eq\begin{cases}(-A^2-4)^{(p-1)/4}\pmod p&\t{if}\ p\eq1\pmod4,
\\(-A^2-4)^{(p+1)/4}u_{(p+1)/2}(A)/2\pmod p&\t{if}\ p\eq3\pmod4.\end{cases}
\end{equation}
\end{theorem}

We will prove Theorem \ref{Th1.4} in Section 4.

Let $p$ be a prime with $p\eq1\pmod4$. Then $(\f{p-1}2!)^2\eq-1\pmod p$ by Wilson's theorem. We may write $p=x^2+y^2$ with $x,y\in\Z$, $x\eq1\pmod4$ and $y\eq \f{p-1}2!x\pmod p$.
Recall that an integer $a$ not divisible by $p$
is a quartic residue modulo $p$ (i.e., $z^4\eq a\pmod p$ for some $z\in\Z$) if and only if $a^{(p-1)/4}\eq1\pmod p$.
 Dirichlet proved that
$2$ is a quartic residue modulo $p$ if and only if $8\mid y$ (see, e.g., \cite[p.\,64, Exer. 28]{IR}).
On the other hand, we have
$$\l|\l\{1\ls k< \f p4:\ \l(\f kp\r)=1\r\}\r|\eq0\pmod2\iff y\eq(-1)^{(p-1)/4}-1\pmod 8$$
as discovered by K. Burde \cite{B} and re-proved by K. S. Williams \cite{W}.
In view of  Williams and J. D. Currie \cite[(1.4)]{WC}, we have
$$2^{(p-1)/4}\eq(-1)^{|\{1\ls k< \f p4:\ (\f kp)=-1\}|}\times\begin{cases}1\pmod p&\t{if}\ p\eq1\pmod8,
\\\f{p-1}2!\pmod p&\t{if}\ p\eq5\pmod 8.\end{cases}$$
By Dirichlet's class number formula (see, e.g., L.E. Dickson \cite[p.\,101]{D}),
$$\f{p-1}2-4\l|\l\{1\ls k<\f p4:\ \l(\f kp\r)=-1\r\}\r|=h(-4p),$$
 where $h(d)$ with $d\eq0,1\pmod4$ not a square denotes the class number of the quadratic field
with discriminant $d$. In 1905, Lerch (see, e.g., \cite{HW}) proved that
$$h(-3p)=2\sum_{1\ls k<p/3}\l(\f kp\r).$$
By \cite[Lemma 14]{WC},
$$(-3)^{(p-1)/4}\eq\begin{cases}(-1)^{h(-3p)/4}\pmod p&\t{if}\ p\eq1\pmod{12},
\\(-1)^{(h(-3p)-2)/4}\f{p-1}2!\pmod p&\t{if}\ p\eq5\pmod {12}.\end{cases}$$
Thus, if $p\eq1\pmod{12}$ then
$$(-3)^{(p-1)/4}\eq(-1)^{\f12\sum_{k=1}^{(p-1)/3}((\f kp)-1)+\f{p-1}6}=(-1)^{|\{1\ls k<\f p3:\ (\f kp)=-1\}|}\pmod p;$$
similarly, if $p\eq5\pmod{12}$ then
$$(-3)^{(p-1)/4}\eq(-1)^{|\{1\ls k<\f p3:\ (\f kp)=-1\}|}\,\f{p-1}2!\pmod p.$$

From Theorem \ref{Th1.4}, we deduce the following result which will be proved in Section 5.

\begin{theorem}\label{Th1.5} Let $p$ be an odd prime.

{\rm (i)} We have
\begin{equation}\label{1.17} T_p(1,-1,-1)
\eq\begin{cases}-5^{(p-1)/4}\pmod p&\t{if}\ p\eq1,9\pmod{20},
\\(-5)^{(p-1)/4}\pmod p&\t{if}\ p\eq 13,17\pmod{20},
\\(-1)^{\lfloor(p-10)/20\rfloor}\pmod p&\t{if}\ p\eq3,7\pmod{20},
\\(-1)^{\lfloor(p-5)/10\rfloor}\pmod p&\t{if}\ p\eq11,19\pmod{20}.
\end{cases}\end{equation}

{\rm (ii)} We have
\begin{equation}\label{1.18} T_p(1,-2,-1)
\eq\begin{cases}-2^{(p-1)/4}\pmod p&\t{if}\ p\eq1\pmod{8},
\\2^{(p-1)/4}\pmod p&\t{if}\ p\eq5\pmod 8,
\\(-1)^{(p-3)/8}\pmod p&\t{if}\ p\eq3\pmod 8,
\\(-1)^{(p-7)/8}\pmod p&\t{if}\ p\eq7\pmod{8}.
\end{cases}\end{equation}
\end{theorem}

Now we state our sixth theorem.

\begin{theorem}\label{Th1.6} Let $p>3$ be a prime and let $\da\in\{\pm1\}$.
If $p\eq1\pmod4$, then
\begin{equation}\label{1.19} T_p(2,5\da,2)\eq(-1)^{\lfloor(p+11)/12\rfloor}\pmod p.
\end{equation}
When $p\eq3\pmod4$, we have
\begin{equation}\label{1.20} T_p(2,5\da,2)
\eq\l(\f 6p\r)\f{\da2^{\da}}{3^{\da}}\bi{(p-3)/2}{(p-3)/4}^{-2\da}\pmod p.
\end{equation}
\end{theorem}

Note that there is no simple closed form for $\bi{(p-3)/2}{(p-3)/4}$ modulo
 a prime $p\eq3\pmod 4$. For a prime $p\eq1\pmod 4$ with $p=x^2+y^2$ ($x,y\in\Z$ and $x\eq1\pmod4$),
Gauss showed the congruence $\bi{(p-1)/2}{(p-1)/4}\eq2x\pmod p$, and  S. Chowla, B. Dwork and R. J. Evans \cite{CDE} used Gauss and Jacobi sums to prove further that
$$\bi{(p-1)/2}{(p-1)/4}\eq\f{2^{p-1}+1}2\l(2x-\f p{2x}\r)\pmod {p^2},$$
which was first conjectured by F. Beukers. (See also \cite[Chapter 9]{BEW} for further related results.)

Though we have made some numerical tests via a computer,
 we are unable to find general patterns for $T_p(a,b,c)$ modulo $p$, where $p$ is an arbitrary odd prime and $a,b,c$ are arbitrary integers.

Let $p$ be an odd prime. It is known (cf. \cite[(1.6) and (1.7)]{S19b}) that
$$\prod_{1\ls i<j\ls (p-1)/2\atop p\nmid i^2+j^2}(i^2+j^2)\eq
\begin{cases}(-1)^{\lfloor(p-5)/8\rfloor}\pmod p&\t{if}\ p\eq1\pmod4,\\(-1)^{\lfloor(p+1)/8\rfloor}\pmod p&\t{if}\ p\eq3\pmod4.\end{cases}$$
From this we immediately get
 $$\prod_{1\ls i<j\ls (p-1)/2\atop p\nmid i^2+j^2}\l(\f{i^2+j^2}p\r)
 =\begin{cases}1&\t{if}\ p\eq1\pmod4,\\(-1)^{\lfloor(p+1)/8\rfloor}&\t{if}\ p\eq3\pmod4.
 \end{cases}$$
As the product $\prod_{1\ls i<j\ls (p-1)/2}(i^2-j^2)$ modulo $p$
was determined via \cite[(1.5)]{S19b}, we also know the value of the product
$$\prod_{1\ls i<j\ls (p-1)/2}\l(\f{i^2-j^2}p\r)=\prod_{1\ls i<j\ls (p-1)/2}\l(\f{i-j}p\r)\l(\f{i+j}p\r)$$
Motivated by this, we obtain the following result.

\begin{theorem}\label{Th1.7} Let $p>3$ be a prime and let $\da\in\{\pm1\}$. Then
\begin{equation}\label{1.21}\prod_{1\ls i<j\ls(p-1)/2}\l(\f{j+\da i}p\r)
=\begin{cases}(-1)^{|\{0<k<\f p4:\ (\f kp)=\da\}|}&\t{if}\ p\eq1\pmod4,
\\(-1)^{(p-3)/8}&\t{if}\ p\eq3\pmod8,
\\(-1)^{(p+1)/8+(h(-p)+1)/2}&\t{if}\ p\eq7\pmod 8.
\end{cases}\end{equation}
\end{theorem}

We will prove Theorems \ref{Th1.6} and \ref{Th1.7} in Section 6, and pose ten conjectures in Section 7.

\section{Proofs of Theorems \ref{Th1.1} and \ref{Th1.2}}
\setcounter{theorem}{0}
\setcounter{corollary}{0}
\setcounter{equation}{0}

\medskip
\noindent
{\bf Proof of Theorem \ref{Th1.1}}. For each $k=1,\ldots,(n-1)/2$, we have
$$\l\lfloor\f{kx}n\r\rfloor-\l\lfloor\f{k(x-1)}n\r\rfloor=\begin{cases}0&\t{if}\ \{kx\}_n>k,
\\1&\t{if}\ \{kx\}_n<k.\end{cases}$$
Thus
$$\l|\l\{1\ls k<\f n2:\ \{kx\}_n>k\r\}\r|
=\f{n-1}2-\sum_{k=1}^{(n-1)/2}\l(\l\lfloor\f{kx}n\r\rfloor-\l\lfloor\f{k(x-1)}n\r\rfloor\r)$$
and hence
\begin{align*}&(-1)^{|\{1\ls k<n/2:\ \{kx\}_n>k\}|}\l(\f{-1}n\r)
\\=&(-1)^{\sum_{k=1}^{(n-1)/2}\lfloor kx/n\rfloor+\sum_{k=1}^{(n-1)/2}\lfloor k(x-1)/n\rfloor}
\\=&(-1)^{(\sum_{k=1}^{(n-1)/2}(2x-1)k-\sum_{k=1}^{(n-1)/2}\{kx\}_n-\sum_{k=1}^{(n-1)/2}\{k(x-1)\}_n)/n}
\\=&(-1)^{(n^2-1)/8}(-1)^{\sum_{k=1}^{(n-1)/2}\{kx\}_n+\sum_{k=1}^{(n-1)/2}\{k(x-1)\}_n}.
\end{align*}
As $\{kx\}_n\eq 1+n-\{kx\}_n\pmod 2$ for all $k=1,\ldots,(n-1)/2$, we have
\begin{align*} \sum_{k=1}^{(n-1)/2}\{kx\}_n\eq&\sum_{k=1\atop \{kx\}_n<n/2}^{(n-1)/2}\{kx\}_n
+\sum_{k=1\atop\{kx\}_n>n/2}^{(n-1)/2}(1+(n-\{kx\}_n))
\\=&\sum_{k=1\atop\{kx\}_n>n/2}^{(n-1)/2}1+\sum_{r=1}^{(n-1)/2}r
\\=&\l|\l\{1\ls k<\f n2:\ \{kx\}_n>\f n2\r\}\r|+\f{n^2-1}8\pmod2.
\end{align*}
and hence
$$(-1)^{\sum_{k=1}^{(n-1)/2}\{kx\}_n}=\l(\f xn\r)\l(\f 2n\r)=\l(\f{2x}n\r)$$
with the help of (1.2). Similarly,
$$(-1)^{\sum_{k=1}^{(n-1)/2}\{k(x-1)\}_n}=\l(\f xn\r)\l(\f 2n\r)=\l(\f{2(x-1)}n\r).$$

In view of the above, we obtain
\begin{align*}(-1)^{|\{1\ls k<n/2:\ \{kx\}_n>k\}|}=&\l(\f{-2}n\r)(-1)^{\sum_{k=1}^{(n-1)/2}\{kx\}_n}(-1)^{\sum_{k=1}^{(n-1)/2}\{k(x-1)\}_n}
\\=&\l(\f{-2}n\r)\l(\f{2x}n\r)\l(\f{2(x-1)}n\r)=\l(\f{2x(1-x)}n\r).
\end{align*}
This proves (\ref{1.3}).

When $1\ls k<n/2$ and $\{kx\}_n<n/2$, we clearly have
$$\{kx\}_n>k\iff\{k(1-x)\}_n>\f n2.$$
Thus
\begin{align*}&\l|\l\{1\ls k<\f n2:\ \{kx\}_n>k\r\}\r|-\l|\l\{1\ls k<\f n2:\ \{kx\}_n>\f n2\r\}\r|
\\=&\l|\l\{1\ls k<\f n2:\ \{kx\}_n<\f n2<\{k(1-x)\}_n\r\}\r|
\\=&\l|\l\{1\ls k<\f n2:\ \{k(1-x)\}_n>\f n2\r\}\r|
\\&-\l|\l\{1\ls k<\f n2:\ \{kx\}_n>\f n2\ \&\ \{k(1-x)\}_n>\f n2\r\}\r|,
\end{align*}
and hence by (\ref{1.2}) and (\ref{1.3}) we get
\begin{align*}&(-1)^{|\{1\ls k<n/2:\ \{kx\}_n>n/2\ \&\ \{k(1-x)\}_n>n/2\}|}
\\=&(-1)^{|\{1\ls k<n/2:\ \{kx\}_n>k\}|}\l(\f x n\r)\l(\f{1-x}n\r)=\l(\f 2n\r).
\end{align*}
This proves (\ref{1.4}).

In view of (\ref{1.2}) and (\ref{1.4}), we also have
\begin{align*}&(-1)^{|\{1\ls k<n/2:\ \{kx\}_n>n/2>\{k(1-x)\}_n\}|}
\\=&(-1)^{|\{1\ls k<n/2:\ \{kx\}_n>n/2\}|-||\{1\ls k<n/2:\ \{kx\}_n>n/2\ \&\ \{k(1-x)\}_n>n/2\}|}
\\=&\l(\f{x}n\r)\l(\f{2}n\r)=\l(\f {2x}n\r).
\end{align*}
This proves (\ref{1.6}).

By (\ref{1.4}) and (\ref{1.6}), we have
\begin{align*} &\l(\f 2n\r)(-1)^{|\{1\ls k<n/2:\ \{kx\}_n<n/2\ \&\ \{k(1-x)\}_n<n/2\}|}
\\=&(-1)^{|\{1\ls k<n/2:\ (\{kx\}_n-n/2)(\{k(1-x)\}_n-n/2)>0\}|}
\\=&(-1)^{(n-1)/2-|\{1\ls k<n/2:\ \{kx\}_n>n/2>\{k(1-x)\}_n\ \t{or}\ \{k(1-x)\}_n>n/2>\{kx\}_n\}|}
\\=&\l(\f{-1}n\r)\l(\f{2x}n\r)\l(\f{2(1-x)}n\r)=\l(\f {x(x-1)}n\r).
\end{align*}
So (\ref{1.5}) also holds.

The proof of Theorem \ref{Th1.1} is now complete. \qed

For any odd prime $p$ and rational $p$-adic integer $x$, we define
\begin{equation}\label{2.1} N_p(x):=|\{1\ls k<p/2:\ \{kx\}_p>k\}|.\end{equation}

\medskip
\noindent{\bf Proof of Theorem \ref{Th1.2}}. By parts (ii)-(iv) of Sun \cite[Theorem 1.2]{S19b},
the desired result holds except for the cases
\begin{align*} \t{I}.&\ p\nmid ac(a-c)\ \t{and}\ p\mid a+b+c,
\\\t{II}.&\ p\nmid ab(a+b)\ \t{and}\ p\mid c,
\\\t{III}.&\ p\mid a\ \t{and}\ p\nmid bc(b+c).
\end{align*}

In case I, by \cite[Theorem 1.2(iii)]{S19b} we have
$$S_p(a,b,c)\eq(-1)^{N_p(a/c)}\l(\f{2c(a-c)}p\r)\pmod p.$$
As
$$(-1)^{N_p(a/c)}=\l(\f{2a(c-a)}p\r)$$
by Theorem \ref{Th1.1}, we obtain
$$S_p(a,b,c)\eq\l(\f{2a(c-a)}p\r)\l(\f{2c(a-c)}p\r)=\l(\f{-ac}p\r)\pmod p.$$

In case II, by \cite[Theorem 1.2(iv)]{S19b} we have
$$S_p(a,b,c)\eq(-1)^{N_p(-a/b)}\l(\f{2}p\r)\pmod p.$$
As
$$(-1)^{N_p(-a/b)}=\l(\f{-2a(a+b)}p\r)$$
by Theorem \ref{Th1.1}, we get
$$S_p(a,b,c)\eq\l(\f{-2a(a+b)}p\r)\l(\f{2}p\r)=\l(\f{-a(a+b)}p\r)\pmod p.$$

In case III, by \cite[Theorem 1.2(iv)]{S19b} we have
$$S_p(a,b,c)\eq(-1)^{N_p(-c/b)}\l(\f{2}p\r)\pmod p.$$
Since $(-1)^{N_p(-c/b)}=(\f{-2c(b+c)}p)$
by Theorem \ref{Th1.1}, we deduce that
$$S_p(a,b,c)\eq\l(\f{-2c(b+c)}p\r)\l(\f{2}p\r)=\l(\f{-c(b+c)}p\r)\pmod p.$$

In view of the above, we have completed the proof of Theorem \ref{Th1.2}. \qed

\section{Proof of Theorem 1.3}
\setcounter{theorem}{0}
\setcounter{corollary}{0}
\setcounter{equation}{0}

\begin{lemma}\label{Lem3.1} Let $p$ be any odd prime. Then
\begin{equation}\label{3.1}\l(\f{p-1}2!\r)^2\eq (-1)^{(p+1)/2}\pmod p,\end{equation}
and
\begin{equation}\label{3.2}\prod_{1\ls i<j\ls(p-1)/2}(j^2-i^2)\eq\begin{cases}-((p-1)/2)!\pmod p&\t{if}\ p\eq1\pmod4,\\1\pmod p&\t{if}\ p\eq3\pmod4.\end{cases}
\end{equation}
\end{lemma}
\begin{remark}\label{Rem3.1} (\ref{3.1}) is an easy consequence of Wilson's theorem. (\ref{3.2}) is also known, see \cite[(1.5)]{S19b} and its few-line proof there.
\end{remark}

\begin{lemma}\label{Lem3.2} Let $p$ be a prime with $p\eq1\pmod 4$. Then
\begin{equation}\label{3.3}\l|\l\{1\ls k\ls\f{p-1}2:\ \l\{k\times\f{p-1}2!\r\}_p>\f p2\r\}\r|=\f{p-1}4.
\end{equation}
\end{lemma}
\Proof. Let $a=((p-1)/2)!$. Then $a^2\eq-1\pmod p$ by (\ref{3.1}). For any $k=1,\ldots,(p-1)/2$, there is a unique integer $k^*\in\{1,\ldots,(p-1)/2\}$ congruent to $ak$ or $-ak$ modulo $p$. Note that $k^*\not=k$ since
$a\not\eq\pm1\pmod p$. If $\{ak\}_p>p/2$ then
$\{ak^*\}_p=\{a(-ak)\}_p=k<p/2$; if $\{ak\}_p<p/2$ then $\{ak^*\}_p=\{a(ak)\}_p=p-k>p/2$.
So, exactly one of $\{ak\}_p$ and $\{ak^*\}_p$ is greater than $p/2$. Therefore (\ref{3.3}) holds. \qed

\begin{remark}\label{Rem3.2} In view of Gauss' Lemma, Lemma \ref{Lem3.2} is stronger than the fact that $(\f{((p-1)/2)!}p)=(\f 2p)$ for any prime $p\eq1\pmod4$ (cf. \cite[Lemma 2.3]{S19a}).
\end{remark}

\medskip
\noindent{\bf Proof of Theorem 1.3}. Observe that
\begin{align*}\prod^{(p-1)/2}_{i,j=1\atop p\nmid i^2-j^2}(i^2-j^2)
=&\prod_{1\ls i<j\ls(p-1)/2}(i^2-j^2)(j^2-i^2)
\\=&(-1)^{\bi{(p-1)/2}2}\prod_{1\ls i<j\ls(p-1)/2}(j^2-i^2)^2
\end{align*}
and hence
\begin{equation}\label{3.4}\prod^{(p-1)/2}_{i,j=1\atop p\nmid i^2-j^2}(i^2-j^2)\eq -\l(\f 2p\r)\pmod p
\end{equation}
with the help of Lemma \ref{Lem3.1}.

When $a\eq\pm1\pmod p$, we have $b\eq\mp1\not\eq a\pmod p$ and hence
$$\{a,b\}_p\eq-\l(\f 2p\r)=\begin{cases}-(\f{a-b}p)\pmod p&\t{if}\ p\eq1\pmod 4,
\\-(\f{a^2+1}p)\pmod p&\t{if}\ p\eq3\pmod4.\end{cases}$$
So (\ref{1.12}) holds in the case $a\eq\pm1\pmod p$.

Below we assume that $a\not\eq \pm1\pmod p$. As $ab\eq-1\pmod p$, we have
\begin{equation}\label{3.5}\begin{aligned} \{a,b\}_p
=&\prod^{(p-1)/2}_{i,j=1\atop i\not\eq aj,bj\pmod p}(i-aj)\times\prod^{(p-1)/2}_{i,j=1\atop j\not\eq ai,bi\pmod p}(j-bi)
\\\eq&\prod^{(p-1)/2}_{i,j=1\atop p\nmid (i-aj)(j+ai)}(i-aj)
\times\prod^{(p-1)/2}_{i,j=1\atop p\nmid (i+aj)(j-ai)}\f{i+aj}a\pmod p.
\end{aligned}\end{equation}

(i) If $p\eq3\pmod 4$, then $(\f{-1}p)=-1$ and hence $a\not\eq b\pmod p$.
Now we prove (\ref{1.12}) under the assumption $a\not\eq b\pmod p$. Note that $a^2\not\eq\pm1\pmod p$.

For $1\ls i,j\ls (p-1)/2$, we cannot have $i^2-a^2j^2\eq j^2-a^2i^2\eq0\pmod p$.
Thus
\begin{align*}&\l|\l\{(i,j):\ 1\ls i,j\ls\f{p-1}2\ \&\  p\nmid(i+aj)(j-ai)\r\}\r|
\\=&\l(\f{p-1}2\r)^2-\l|\l\{(i,j):\ 1\ls i,j\ls\f{p-1}2\ \&\ p\mid i+aj\r\}\r|
\\&-\l|\l\{(i,j):\ 1\ls i,j\ls\f{p-1}2\ \&\ p\mid j-ai\r\}\r|
\\=&\l(\f{p-1}2\r)^2-\l|\l\{1\ls j\ls\f{p-1}2:\ \{aj\}_p>\f p2\r\}\r|
\\&-\l|\l\{1\ls i\ls\f{p-1}2:\ \{ai\}_p<\f p2\r\}\r|
\\=&(p-1)\f{p-1}4-\f{p-1}2=\f{p-1}2\cdot\f{p+1}2-(p-1)
\end{align*}
and hence
$$a^{|\{(i,j):\ 1\ls i,j\ls(p-1)/2\ \&\ p\nmid(i+aj)(j-ai)\}|}\eq (a^{(p-1)/2})^{(p+1)/2}\eq\l(\f ap\r)^{(p+1)/2}\pmod p.$$
Combining this with (\ref{3.5}), we obtain
\begin{align*}&\l(\f ap\r)^{(p+1)/2}\{a,b\}_p
\\\eq&\prod^{(p-1)/2}_{i,j=1\atop p\nmid (i^2-a^2j^2)(j^2-a^2i^2)}(i^2-a^2j^2)
\times \prod^{(p-1)/2}_{i,j=1\atop i\eq aj\pmod p}(i+aj)
\times \prod^{(p-1)/2}_{i,j=1\atop i\eq -aj\pmod p}(i-aj)
\\&\times\prod^{(p-1)/2}_{i,j=1\atop  j-ai\eq-a(i+bj)\eq 0\pmod p}(i-aj)\times\prod^{(p-1)/2}_{i,j=1\atop
j+ai\eq a(i-bj)\eq0\pmod p}(i+aj)
 \\\eq&\prod^{(p-1)/2}_{i,j=1\atop p\nmid i^2-a^2j^2}(i^2-a^2j^2)\bigg/\prod^{(p-1)/2}_{i,j=1\atop p\mid j^2-a^2i^2}(i^2-a^2j^2)
\\&\times(-1)^{|\{1\ls j\ls(p-1)/2:\ \{aj\}_p>p/2\}|}\prod_{j=1}^{(p-1)/2}(2aj)
\\&\times(-1)^{|\{1\ls j\ls(p-1)/2:\ \{bj\}_p>p/2\}|}\prod_{j=1}^{(p-1)/2}(bj+aj)
\pmod p
\end{align*}
Thus, by using (\ref{3.4}) and Gauss' Lemma, we see that
\begin{align*}\l(\f ap\r)^{(p+1)/2}\{a,b\}_p
\eq& \prod^{(p-1)/2}_{i,k=1\atop p\nmid i^2-k^2}(i^2-k^2)\bigg/\prod^{(p-1)/2}_{i=1}(i^2-a^2(a^2i^2))
\\&\times \l(\f{ab}p\r)(2a(a+b))^{(p-1)/2}\l(\f{p-1}2!\r)^2
\\\eq&-\l(\f 2p\r)\l(\f{-1}p\r)\l(\f{2(a^2-1)(1-a^4)}p\r)
\\=&-\l(\f{a^2+1}p\r)=-\l(\f{a^2-ab}p\r)\pmod p.
\end{align*}
This proves (\ref{1.12}).

 (ii) Now we come to prove part (ii) of Theorem \ref{Th1.3}. Assume that $a\eq b\pmod p$.
 Then $a^2\eq ab\eq-1\pmod p$ and hence $p\eq1\pmod4$. As $j\pm ai\eq \pm a(i\mp aj)\pmod p$,
by (\ref{3.5}) we have
\begin{align*} \{a,b\}_p\eq& a^{-|\{(i,j):\ 1\ls i,j\ls(p-1)/2\ \&\ p\nmid i+aj\}|}
\prod^{(p-1)/2}_{i,j=1\atop p\nmid i-aj}(i-aj)\times\prod^{(p-1)/2}_{i,j=1\atop p\nmid i+aj}(i+aj)
\\\eq&a^{-(p-1)^2/4+|\{(i,j):\ 1\ls i,j\ls(p-1)/2\ \&\ p\mid i+aj\}|}
\prod^{(p-1)/2}_{i,j=1\atop p\nmid i^2-(aj)^2}(i^2-(aj)^2)
\\&\times\prod^{(p-1)/2}_{i,j=1\atop p\mid i-aj}(i+aj)\times\prod^{(p-1)/2}_{i,j=1\atop p\mid i+aj}(i-aj)
\\\eq&a^{|\{(i,j):\ 1\ls i,j\ls(p-1)/2\ \&\ p\mid i+aj\}|}\prod^{(p-1)/2}_{i,k=1\atop p\nmid i^2-k^2}(i^2-k^2)
\\&\times(-1)^{|\{1\ls j\ls(p-1)/2:\ \{aj\}_p>p/2\}|}\prod_{j=1}^{(p-1)/2}(2aj).
\pmod p\end{align*}
Applying (\ref{3.4}) and Gauss' Lemma, from the above we get
\begin{equation}\label{3.6}\begin{aligned} \{a,b\}_p\eq&a^{|\{1\ls j\ls(p-1)/2:\ \{aj\}_p>p/2\}|}\times\l(-\l(\f 2p\r)\r)\l(\f ap\r)(2a)^{(p-1)/2}\f{p-1}2!
\\\eq&-\f{p-1}2!\times a^{|\{1\ls j\ls(p-1)/2:\ \{aj\}_p>p/2\}|}\pmod p.
\end{aligned}\end{equation}
As $a^2\eq-1\eq((p-1)/2)!)^2$, we have $a\eq(-1)^k((p-1)/2)!\pmod p$ for some $k\in\{0,1\}$.
In view of Lemma \ref{Lem3.2},
\begin{align*} &\l|\l\{1\ls j\ls \f{p-1}2:\ \l\{-\f{p-1}2!\,j\r\}_p>\f p2\r\}\r|
\\=&\l|\l\{1\ls j\ls \f{p-1}2:\ \l\{\f{p-1}2!\,j\r\}_p<\f p2\r\}\r|=\f{p-1}4.
\end{align*}
Hence
\begin{align*} &a^{|\{1\ls j\ls(p-1)/2:\ \{aj\}_p>p/2\}|}
\\=&a^{(p-1)/4}=(-1)^{k(p-1)/4}\l(\f{p-1}2!\r)^{(p-1)/4}
\\\eq&\begin{cases}(\f{p-1}2!)^{2(p-1)/8}\eq(-1)^{(p-1)/8}\pmod p&\t{if}\ p\eq1\pmod 8,
\\(-1)^k\f{p-1}2!(\f{p-1}2!)^{2(p-5)/8}\eq(-1)^{k+(p-5)/8}\f{p-1}2!\pmod p&\t{if}\ p\eq5\pmod 8.
\end{cases}
\end{align*}
Combining this with (\ref{3.6}) we immediately obtain the desired results in Theorem \ref{Th1.3}(ii).

In view of the above, we have completed the proof of Theorem \ref{Th1.3}. \qed

\section{Proof of Theorem \ref{Th1.4}}
\setcounter{theorem}{0}
\setcounter{corollary}{0}
\setcounter{equation}{0}

\noindent{\bf Proof of Theorem \ref{Th1.4}(i)}. As $A^2\eq-4\pmod p$, we have $(\f{-1}p)=1$ and hence $p\eq1\pmod4$. Since $((p-1)/2)^2\eq-1\eq(A/2)^2\pmod p$, for some $k\in\{0,1\}$ we have $A/2\eq(-1)^k((p-1)/2)!\pmod p$. Choose $a,b\in\Z$ with $a\eq b\eq A/2\pmod p$.
Note that $\{a,b\}_p\eq T_p(1,-A,-1)\pmod p$.
 Applying Theorem \ref{Th1.3}(ii) we immediately get the desired (\ref{1.14}). \qed
 \medskip

For any odd prime $p$ and integer $A$ with $\Delta=A^2+4\not\eq0\pmod p$, it is known (cf. \cite[Lemma 2.3]{S10}) that
$$u_{(p-(\f{\Delta}p)/2}(A)v_{(p-(\f{\Delta}p))/2}(A)=u_{p-(\f{\Delta}p)}(A)\eq0\pmod p.$$

\begin{lemma}\label{Lem4.1} Let $A\in\Z$ and let $p$ be an odd prime not dividing $\Delta=A^2+4$. Then
\begin{equation}\label{4.1}p\mid v_{(p-(\f{\Delta}p))/2}(A)\iff \l(\f{-1}p\r)=-1.
\end{equation}
\end{lemma}
\begin{remark}\label{Rem3.1} This is a known result, see, e.g., \cite[Chapter 2, (IV.23)]{Ri}.
\end{remark}

\medskip
\noindent{\bf Proof of Theorem \ref{Th1.4}(ii)}. Let $\Delta=A^2+4$. As $(\frac{\Delta}p)=1$, we have
$d^2\eq \Delta$ for some $d\in\Z$ with $p\nmid d$.
Choose integers $a$ and $b$ such that
$$a\eq\f{A+d}2\pmod p\ \t{and}\ b\eq\f{A-d}2\pmod p.$$
Then $a+b\eq A\pmod p$ and $ab\eq(A^2-d^2)/4\eq-1\pmod p$. Thus, for any $i,j\in\Z$ we have
$$i^2-Aij-j^2\eq(i-aj)(i-bj)\pmod p.$$
If $p\eq1\pmod4$, then
$$(A^2+4)^{(p-1)/4}\eq d^{(p-1)/2}\eq(a-b)^{(p-1)/2}\eq\l(\f{a-b}p\r)\pmod p.$$
If $p\eq3\pmod4$, then
$$\l(\f{2d(A+d)}p\r)=\l(\f{ad}p\r)=\l(\f{a(a-b)}p\r).$$

For any positive integer $n$, clearly
\begin{align*} u_n(A)=&\f1{\sqrt{\Delta}}\(\l(\f{A+\sqrt{\Delta}}2\r)^n-\l(\f{A-\sqrt{\Delta}}2\r)^n\)
\\=&\f1{2^{n-1}}\sum_{k=0}^{\lfloor(n-1)/2\rfloor}\bi n{2k+1}A^{n-1-2k}\Delta^k
\\\eq&\f1{2^{n-1}}\sum_{k=0}^{\lfloor(n-1)/2\rfloor}\bi n{2k+1}A^{n-1-2k}d^{2k}
\\=&\f1{d}\(\l(\f{A+d}2\r)^n-\l(\f{A-d}2\r)^n\)\pmod p
\end{align*}
and
\begin{align*} v_n(A)=&\l(\f{A+\sqrt{\Delta}}2\r)^n+\l(\f{A-\sqrt{\Delta}}2\r)^n
=\f1{2^{n-1}}\sum_{k=0}^{\lfloor n/2\rfloor}\bi n{2k}A^{n-2k}\Delta^k
\\\eq&\f1{2^{n-1}}\sum_{k=0}^{\lfloor n/2\rfloor}\bi n{2k}A^{n-2k}d^{2k}
=\l(\f{A+d}2\r)^n+\l(\f{A-d}2\r)^n\pmod p.
\end{align*}
In the case $p\eq3\pmod4$, we have $p\mid v_{(p-1)/2}(A)$ by Lemma \ref{Lem4.1}, and hence
\begin{align*}\l(\f{2\da(A+d)}p\r)\eq&d^{(p-1)/2}\l(\f{A+d}2\r)^{(p-1)/2}
\eq d^{(p-1)/2}\f{v_{(p-1)/2}(A)+d u_{(p-1)/2}(A)}2
\\\eq&\l(d^2\r)^{(p+1)/4}\f{u_{(p-1)/2}(A)}2
\eq\Delta^{(p+1)/4}\f{u_{(p-1)/2}(A)}2\pmod p.
\end{align*}

In view of the above, we obtain (\ref{1.15}) from Theorem \ref{Th1.3}(i). \qed

\begin{lemma}\label{Lem4.2} Let $p$ be an odd prime, and let $A\in\Z$ with $\Delta=A^2+4\not\eq0\pmod p$. If $p\eq1\pmod4$, then
\begin{equation}\label{4.2}u_{(p+(\f{\Delta}p))/2}(A)\eq \pm\Delta^{(p-1)/4}\pmod p.
\end{equation}
If $p\eq3\pmod4$, then
\begin{equation}\label{4.2}u_{(p-(\f{\Delta}p))/2}(A)\eq \pm2\Delta^{(p-3)/4}\pmod p.
\end{equation}
\end{lemma}
\begin{remark}\label{Rem4.2} This is a known result, see, e.g., \cite[Theorem 4.1]{Su}.
\end{remark}

\medskip\noindent{\bf Proof of Theorem \ref{Th1.4}(iii)}.
Let
 $$\Delta= A^2+4,\ \al=\f{A+\sqrt{\Delta}}2\ \t{and}\ \beta=\f{A-\sqrt{\Delta}}2.$$
 As $x^2-Ax-1=(x-\al)(x-\beta)$, both $\al$ and $\beta$ are algebraic integers.
 Observe that
 \begin{align*}&\prod_{i,j=1}^{(p-1)/2}(i^2-Aij-j^2)
 \\=&\prod_{i,j=1}^{(p-1)/2}(i-\al j)(i-\beta j)
 =\prod_{i,j=1}^{(p-1)/2}(i-\al j)\times\prod_{i,j=1}^{(p-1)/2}(j-\beta i)
 \\=&\prod_{i,j=1}^{(p-1)/2}(i-\al j)\times\prod_{i,j=1}^{(p-1)/2}\f{\al j+i}{\al}
 =\al^{-((p-1)/2)^2}\prod_{i,j=1}^{(p-1)/2}(i^2-\al^2j^2)
 \\=&\al^{-((p-1)/2)^2}\prod_{i,j=1}^{(p-1)/2}(-j^2)\l(\al^2-\f{i^2}{j^2}\r)
 \\=&(-\al)^{-((p-1)/2)^2}\l(\f{p-1}2!\r)^{2\f{p-1}2}\prod_{j=1}^{(p-1)/2}
 \prod_{i=1}^{(p-1)/2}\l(\al^2-\f{i^2}{j^2}\r)
 \end{align*}
 and hence
 $$\prod_{i,j=1}^{(p-1)/2}(i^2-Aij-j^2)
 \eq\beta^{((p-1)/2)^2}\prod_{k=1}^{(p-1)/2}(\al^2-k^2)^{(p-1)/2}
\pmod p$$
(in the ring of all algebraic integers) since for any $j,k=1,\ldots,(p-1)/2$ there is a unique $i\in\{1,\ldots,(p-1)/2\}$ congruent to $jk$ or $-jk$ modulo $p$. Note that
$$\prod_{x=1}^{(p-1)/2}(x-k^2)\eq x^{(p-1)/2}-1\pmod p.$$
Thus
\begin{align*}\prod_{i,j=1}^{(p-1)/2}(i^2-Aij-j^2)
 \eq&\beta^{((p-1)/2)^2}((\al^2)^{(p-1)/2}-1)^{(p-1)/2}
 \\=&\beta^{((p-1)/2)^2}\l(\al^{p-1}-1\r)^{(p-1)/2}
 \\=&\l((-\al)^{(p-1)/2}-\beta^{(p-1)/2}\r)^{(p-1)/2}
 \pmod p.\end{align*}

 {\bf Case 1.} $p\eq1\pmod4$.

 In this case,
 \begin{align*} \prod_{i,j=1}^{(p-1)/2}(i^2-Aij-j^2)
 \eq&\l(\f{\al^{(p-1)/2}-\beta^{(p-1)/2}}{\al-\beta}\r)^{(p-1)/2}(\al-\beta)^{(p-1)/2}
 \\=&u_{(p-1)/2}(A)^{(p-1)/2}\Delta^{(p-1)/4}\pmod p.
 \end{align*}
 As $(\f{\Delta}p)=-1$, by Lemma \ref{Lem4.2} we have
 $u_{(p-1)/2}(A)\eq\pm\Delta^{(p-1)/4}\pmod p$ and hence
 $$u_{(p-1)/2}(A)^{(p-1)/2}\eq\Delta^{\f{p-1}2\cdot\f{p-1}4}\eq(-1)^{(p-1)/4}\pmod p.$$
 Therefore
 $$\prod_{i,j=1}^{(p-1)/2}(i^2-Aij-j^2)\eq(-\Delta)^{(p-1)/4}\pmod p.$$

  {\bf Case 2.} $p\eq3\pmod4$.

 In this case,
 \begin{align*}\prod_{i,j=1}^{(p-1)/2}(i^2-Aij-j^2)
 \eq&-\l(\al^{(p-1)/2}+\beta^{(p-1)/2}\r)^{(p-1)/2}
 \\=&-v_{(p-1)/2}(A)^{(p-1)/2}\pmod p.
 \end{align*}
 It is easy to see that
 $2v_{n-1}(A)=\Delta u_n(A)-Av_n(A)$ for all $n=1,2,3,\ldots$.
 Hence
 $$2v_{(p-1)/2}(A)=\Delta u_{(p+1)/2}(A)-Av_{(p+1)/2}(A)\eq\Delta u_{(p+1)/2}(A)\pmod p$$
 since $v_{(p+1)/2}(A)=v_{(p-(\f{\Delta}p))/2}(A)\eq0\pmod p$ by Lemma \ref{Lem4.1}.
 Thus
 \begin{align*}\prod_{i,j=1}^{(p-1)/2}(i^2-Aij-j^2)
 \eq&-\l(\f{\Delta}2u_{(p+1)/2}(A)\r)^{(p-1)/2}
 \\\eq&\f{u_{(p+1)/2}(A)}2\l(\f{u_{(p+1)/2}(A)}2\r)^{(p-3)/2}\pmod p.
 \end{align*}
 As $(\f{\Delta}p)=-1$, by Lemma \ref{Lem4.2} we have
 $u_{(p+1)/2}(A)\eq\pm2\Delta^{(p-3)/4}\pmod p$, and hence
 \begin{align*}\prod_{i,j=1}^{(p-1)/2}(i^2-Aij-j^2)
 \eq&\f{u_{(p+1)/2}(A)}2\Delta^{\f{p-3}4\cdot\f{p-3}2}
 =\f{u_{(p+1)/2}(A)}2\Delta^{\f{p+1}2\cdot\f{p+1}4-(p-1)}
 \\\eq&\f{u_{(p+1)/2}(A)}2(-\Delta)^{(p+1)/4}\pmod p.
 \end{align*}

 In view of the above, we have completed the proof of Theorem \ref{Th1.4}(iii). \qed

\section{Proof of Theorem \ref{Th1.5}}
\setcounter{theorem}{0}
\setcounter{corollary}{0}
\setcounter{equation}{0}

\medskip
\noindent{\bf Proof Theorem 1.5}. Note that $\{u_n(1)\}_{n\gs0}$ is just the Fibonacci sequence $\{F_n\}_{n\gs0}$. By \cite[Corollary 2(iii)]{SS}, if $p\eq3\pmod4$ and $(\f 5p)=1$, then
$$u_{(p-1)/2}(1)=F_{(p-1)/2}\eq-2(-1)^{\lfloor(p-5)/10\rfloor}5^{(p-3)/4}\pmod p$$
and hence
$$-5^{(p+1)/4}\f{u_{(p-1)/2}(1)}2\eq(-1)^{\lfloor(p-5)/10\rfloor}5^{(p+1)/4+(p-3)/4}
\eq(-1)^{\lfloor(p-5)/10\rfloor}
\pmod p.$$
Similarly, if $p\eq3\pmod4$ and $(\f 5p)=-1$, then
$$u_{(p+1)/2}(1)=F_{(p+1)/2}\eq2(-1)^{\lfloor(p-5)/10\rfloor}5^{(p-3)/4}\pmod p$$
and hence
$$5^{(p+1)/4}\f{u_{(p+1)/2}(1)}2\eq(-1)^{\lfloor(p-5)/10\rfloor}5^{(p+1)/4+(p-3)/4}
\eq(-1)^{\lfloor (p+5)/10\rfloor}
\pmod p.$$
Therefore Theorem \ref{Th1.4} with $A=1$ yields (\ref{1.17}).

When $p\eq1\pmod 4$, we obviously have
$$8^{(p-1)/4}=2^{(p-1)/2+(p-1)/4}\eq\l(\f 2p\r)2^{(p-1)/4}=(-2)^{(p-1)/4}\pmod p.$$
If $p\eq7\pmod8$, then
$$u_{(p-1)/2}(2)\eq(-1)^{(p+1)/8}2^{(p-3)/4}\pmod p$$
by \cite[(1.7)]{Su},
and hence
\begin{align*}-8^{(p+1)/4}\f{u_{(p-1)/2}(2)}2\eq&-2^{(3p+3)/4}(-1)^{(p+1)/8}2^{(p-7)/4}
\eq(-1)^{(p-7)/8}\pmod p.\end{align*}
Similarly, if $p\eq3\pmod8$, then
$$u_{(p+1)/2}(2)\eq(-1)^{(p+5)/8}2^{(p-3)/4}\pmod p$$
by \cite[(1.7)]{Su},
and hence
\begin{align*}(-8)^{(p+1)/4}\f{u_{(p+1)/2}(2)}2\eq&-2^{(3p+3)/4}(-1)^{(p+5)/8}2^{(p-7)/4}
\eq(-1)^{(p-3)/8}\pmod p.\end{align*}
So Theorem \ref{Th1.4} with $A=2$ yields (1.18). This ends the proof. \qed

\section{Proofs of Theorems \ref{Th1.6} and \ref{Th1.7}}
\setcounter{theorem}{0}
\setcounter{corollary}{0}
\setcounter{equation}{0}

For $n=1,2,3,\ldots$, we adopt the notation
$$n!!:=\prod_{k=0}^{\lfloor (n-1)/2\rfloor}(n-2k).$$

\begin{lemma}\label{Lem6.1} Let $p>3$ be a prime. Then
\begin{equation}\label{6.1}\f{p-1}2!!\prod^{(p-1)/2}_{i,j=1\atop p\nmid 2i+j}(2i+j)
\eq\l(\f {-2}p\r)\f{p-3}2!!\prod^{(p-1)/2}_{i,j=1\atop p\nmid 2i-j}(2i-j)\eq\pm1\pmod p.
\end{equation}
\end{lemma}
\Proof. Set
$$A_p:=\f{p-1}2!!\prod^{(p-1)/2}_{i,j=1\atop p\nmid 2i+j}(2i+j)
\ \ \t{and}\ \ B_p:=\f{p-3}2!!\prod^{(p-1)/2}_{i,j=1\atop p\nmid 2i-j}(2i-j).$$
Then
\begin{align*} \f{A_pB_p}{((p-1)/2)!}\eq&\prod^{(p-1)/2}_{i,j=1\atop p\nmid 2i+j}(2i+j)\times\prod^{(p-1)/2}_{i,j=1\atop p\nmid 2i+p-j}(2i+p-j)
=\prod_{i=1}^{(p-1)/2}\(\f1{2i}\prod^{p-1}_{j=0\atop p\nmid 2i+j}(2i+j)\)
\\\eq&\prod_{i=1}^{(p-1)/2}\f{(p-1)!}{2i}\eq\f{(\f{-2}p)}{((p-1)/2)!}\pmod p
\end{align*}
and hence
\begin{equation}\label{6.2} A_pB_p\eq\l(\f{-2}p\r)\pmod p.\end{equation}

On the other hand,
\begin{align*}\f{B_p}{((p-3)/2)!!}=&\prod^{(p-1)/2}_{i,j=1\atop p\nmid 2(\f{p+1}2-i)-j}\l(2\l(\f{p+1}2-i\r)-j\r)
\\\eq&\prod^{(p-1)/2}_{i,j=1\atop p\nmid 2i+j-1}(1-j-2i)=\prod_{i=1}^{(p-1)/2}\f{-2i}{(-(p-1)/2-2i)^*}\prod^{(p-1)/2}_{j=1\atop p\nmid j+2i}(-j-2i)
\\\eq&\f{(-2)^{(p-1)/2}((p-1)/2)!}{\prod_{1\ls i<p/2\atop 4i\not=p+1}((p+1)/2-2i)}\prod^{(p-1)/2}_{i,j=1\atop p\nmid 2i+j}(2i+j)
\\&\times(-1)^{\sum_{i=1}^{(p-1)/2}((p-1)/2-|\{1\ls j\ls(p-1)/2:\ p\mid 2i+j\}|)}
\pmod p,
\end{align*}
where $k^*$ is $1$ or $k$ according as $p\mid k$ or not.
Note that
$$\prod_{1\ls i<p/4}\l(\f{p+1}2-2i\r)=\f{p-3}2!!.$$
Therefore
\begin{align*}\l(\f{-2}p\r)\f{B_p}{A_p}\eq&\f{((p-1)/2)!}{((p-1)/2)!!}\times\prod_{(p+1)/4<i<p/2}\f1{((p+1)/2-2i)}
\\&\times(-1)^{((p-1)/2)^2-|\{1\ls i\ls(p-1)/2:\ 2i>p/2\}|}
\\\eq&\f{((p-3)/2)!!(-1)^{((p-1)/2)^2-((p-1)/2-\lfloor(p-1)/4\rfloor)}}{(-1)^{\lfloor(p-1)/4\rfloor}\prod_{(p+1)/4<i<p/2}(2i-(p+1)/2)}=1
\pmod p,
\end{align*}
and hence
\begin{equation}\label{6.3} A_p\eq\l(\f{-2}p\r)B_p\pmod p\end{equation}
which gives the first congruence in (\ref{6.1}).

Combining (\ref{6.2}) and (\ref{6.3}), we see that $A_p^2\eq B_p^2\eq1\pmod p$.
So (\ref{6.1}) does hold. This ends the proof. \qed

\medskip
\noindent{\bf Proof of Theorem \ref{Th1.6}}.
As $2i^2+\da5ij+2j^2=(i+\da 2j)(2i+\da j)$, we have
\begin{align*}\prod^{(p-1)/2}_{i,j=1\atop p\nmid 2i^2+\da 5ij+2j^2}(2i^2+\da 5ij+2j^2)
=&\prod^{(p-1)/2}_{i,j=1\atop p\nmid(i+\da2j)(2i+\da j)}(i+\da2j)
\times\prod^{(p-1)/2}_{i,j=1\atop p\nmid(i+\da2j)(2i+\da j)}(2i+\da j)
\\=&\prod^{(p-1)/2}_{i,j=1\atop p\nmid(i+\da2j)(2i+\da j)}(i+\da 2j)
\times\prod^{(p-1)/2}_{i,j=1\atop p\nmid(i+\da 2j)(2i+\da j)}(2j+\da i)
\\=&\prod^{(p-1)/2}_{i,j=1\atop p\nmid i+\da 2j}\da(i+\da 2j)^2\times\prod_{i,j=1\atop p\mid j+\da 2i}^{(p-1)/2}\f1{\da(i+\da 2j)^2}.\end{align*}
(Note that $i+\da 2j\eq j+\da 2i\eq0\pmod p$ for no $i,j=1,\ldots,(p-1)/2$.)
Thus, applying Lemma \ref{Lem6.1} and (\ref{3.1}) we get
\begin{align*}&\prod^{(p-1)/2}_{i,j=1\atop p\nmid 2i^2+\da5ij+2j^2}(2i^2+\da5ij+2j^2)
\\\eq&\f{{\da}^{|(i,j):\ 1\ls i,j\ls (p-1)/2\ \&\ p\nmid i+\da2j\}|}}{((p-2+\da)/2)!!^2}\times\prod_{i=1\atop\{\da2i\}_p>p/2}^{(p-1)/2}\f1{\da(i+\da2(p-2\da i))^2}
\\\eq&\f{{\da}^{((p-1)/2)^2-|\{1\ls j\ls (p-1)/2:\ \{\da2j\}_p>p/2\}|}}
{((p-2+\da)/2)!!^2}\times\prod_{i=1\atop\{\da2i\}_p>p/2}^{(p-1)/2}\f1{\da(i-4i)^2}
\\\eq&\f{\da^{(p-1)^2/4}}{((p-2+\da)/2)!!^2}\times\prod_{i=1}^{(p-1)/2}(3i)^{-1-\da}
\times\prod_{i=1}^{\lfloor p/4\rfloor}(3i)^{2\da}
\\\eq&\f{\da^{(p-1)^2/4}}{((p-2+\da)/2)!!^2}\times\f{3^{2\da\lfloor p/4\rfloor}\lfloor p/4\rfloor!^{2\da}}{((p-1)/2)!^{1+\da}}
\\\eq&\da^{(p-1)^2/4}(-1)^{\f{p+1}2\cdot\f{1+\da}2}3^{2\da\lfloor p/4\rfloor}\l(\f{\lfloor p/4\rfloor!^{\da}}{((p-2+\da)/2)!!}\r)^2
\\=&(-1)^{(p+\da)/2}3^{2\da\lfloor p/4\rfloor}\l(\f{\lfloor p/4\rfloor!^{\da}}{((p-2+\da)/2)!!}\r)^2
\pmod p.
\end{align*}

{\bf Case 1.} $p\eq1\pmod4$.

In this case,
$$\l(\f{\lfloor p/4\rfloor!^{\da}}{((p-2+\da)/2)!!}\r)^2=\begin{cases} 1/2^{(p-1)/2}&\t{if}\ \da=1,
\\2^{(p-1)/2}/((p-1)/2)!^2&\t{if}\ \da=-1.\end{cases}$$
Thus, with the help of (\ref{3.1}) we have
\begin{align*}&(-1)^{(p+\da)/2}3^{2\da\lfloor p/4\rfloor}\l(\f{\lfloor p/4\rfloor!^{\da}}{((p-2+\da)/2)!!}\r)^2
\\\eq&(-1)^{(1+\da)/2}\l(\f 3p\r)^{\da}\l(\f 2p\r)^{-\da}\da=-\l(\f 6p\r)=(-1)^{\lfloor (p+11)/12\rfloor}\pmod p,
\end{align*}
and hence (\ref{1.19}) follows from the above.

{\bf Case 2.} $p\eq3\pmod4$.

In this case, with the aid of (\ref{3.1}) we have
\begin{align*}&(-1)^{(p+\da)/2}3^{2\da\lfloor p/4\rfloor}\l(\f{\lfloor p/4\rfloor!^{\da}}{((p-2+\da)/2)!!}\r)^2
\\\eq&(-1)^{(\da-1)/2}3^{\da((p-1)/2-1)}\l(\f{2^{(p-3)/4}}{(p-1)/2\times\bi{(p-3)/2}{(p-3)/4}}\r)^{2\da}
\l(\f{p-1}2!\r)^{\da-1}
\\\eq&\da\l(\f 3p\r)3^{-\da}2^{\da(p+1)/2}\bi{(p-3)/2}{(p-3)/4}^{-2\da}
\eq\l(\f 6p\r)\f{\da2^{\da}}{3^{\da}}\bi{(p-3)/2}{(p-3)/4}^{-2\da}\pmod p,
\end{align*}
and hence (\ref{1.20}) holds.

In view of the above, we have completed the proof. \qed

\medskip
\noindent{\bf Proof of Theorem \ref{Th1.7}}. Let $n=(p-1)/2$. Clearly,
$$\prod_{1\ls i<j\ls n}(j-i)=\prod_{k=1}^nk^{|\{1\ls i\ls n:\ i+k\ls n\}|}=\prod_{k=1}^nk^{n-k}$$
and hence
\begin{equation}\label{6.4}\prod_{1\ls i<j\ls n}\l(\f{j-i}p\r)=\l(\f{n!}p\r)^n\prod_{k=1\atop 2\nmid k}^n\l(\f kp\r).\end{equation}

{\bf Case 1.} $p\eq1\pmod4$.

In this case, $n$ is even, and hence by (6.4) and \cite[(3.5)]{S19b} we have
$$\prod_{1\ls i<j\ls n}\l(\f{j-i}p\r)=\l(\f{(n-1)!!}p\r)=(-1)^{|\{0<k<\f p4:\ (\f kp)=-1\}|}.$$
By (\ref{3.2}) and \cite[Lemma 2.3]{S19a},
$$\prod_{1\ls i<j\ls n}\l(\f{j^2-i^2}p\r)=\l(\f{-n!}p\r)=\l(\f 2p\r)=(-1)^{(p-1)/4}.$$
Thus we also have
$$\prod_{1\ls i<j\ls n}\l(\f{j+i}p\r)=(-1)^{(p-1)/4}(-1)^{|\{0<k<\f p4:\ (\f kp)=-1\}|}
=(-1)^{|\{0<k<\f p4:\ (\f kp)=1\}|}.$$
So (\ref{1.21}) holds in this case.

{\bf Case 2.} $p\eq3\pmod4$.

In this case, $n$ is odd and
$$\prod_{1\ls i<j\ls n}\l(\f{j^2-i^2}p\r)=\l(\f1p\r)=1$$
by (\ref{3.2}). So it suffices to prove (\ref{1.21}) for $\da=-1$.

In view of (\ref{6.4}), we have
\begin{equation}\label{6.5} \prod_{1\ls i<j\ls n}\l(\f{j-i}p\r)=\l(\f{n!}p\r)\l(\f{n!!}p\r)=\l(\f{(n-1)!!}p\r).
\end{equation}
If $p\eq3\pmod8$, then by \cite[(3.6)]{S19b} we have
$$\l(\f{(n-1)!!}p\r)=(-1)^{\lfloor (p+1)/8\rfloor}=(-1)^{(p-3)/8}$$
and hence (\ref{1.21}) holds for $\da=-1$. When $p\eq7\pmod 8$, we have
$$\l(\f{n!}p\r)=\l(\f{(-1)^{(h(-p)+1)/2}}p\r)=(-1)^{(h(-p)+1)/2}
\ \t{and}\ \l(\f{n!!}p\r)=(-1)^{(p+1)/8}$$
by Mordell \cite{M} and Sun \cite[(3.6)]{S19b} respectively, hence (\ref{1.21}) with $\da=-1$ follows from (\ref{6.5}).

Combining the above, we have finished the proof of (\ref{1.21}). \qed

\section{Some related conjectures}
\setcounter{theorem}{0}
\setcounter{corollary}{0}
\setcounter{equation}{0}
\setcounter{conjecture}{0}

Motivated by our results in Section 1, here we pose 10 conjectures for further research.
We have verified all the following conjectures for primes $p<13000$.

\begin{conjecture}\label{Conj7.1} Let $p>3$ be a prime and let $\da\in\{\pm1\}$. Then
\begin{equation}\label{7.1}\prod_{1\ls i<j\ls(p-1)/2\atop p\nmid 2i^2+\da5ij+2j^2}\l(\f{2i^2+\da5ij+2j^2}p\r)
=\f12\l(\f{\da}p\r)\l(\l(\f{-1}p\r)+\l(\f 2p\r)+\l(\f 6p\r)-\l(\f p3\r)\r).
\end{equation}
\end{conjecture}

\begin{conjecture}\label{Conj7.2} Let $p>3$ be a prime. Then
\begin{equation}\label{7.2}\prod_{1\ls i<j\ls(p-1)/2\atop p\nmid i^2-ij+j^2}\l(\f{i^2-ij+j^2}p\r)
=\begin{cases}-1&\t{if}\ p\eq 5,7\pmod {24},\\1&\t{otherwise}.\end{cases}
\end{equation}
Also,
\begin{equation}\label{7.3}\prod_{1\ls i<j\ls(p-1)/2\atop p\nmid i^2+ij+j^2}\l(\f{i^2+ij+j^2}p\r)
=\begin{cases}-1&\t{if}\ p\eq 5,11\pmod {24},\\1&\t{otherwise}.\end{cases}
\end{equation}
\end{conjecture}

\begin{conjecture}\label{Conj7.3} Let $p>3$ be a prime. Then
\begin{equation}\label{7.4}\prod_{1\ls i<j\ls(p-1)/2\atop p\nmid i^2-3ij+j^2}\l(\f{i^2-3ij+j^2}p\r)
=\begin{cases}-1&\t{if}\ p\eq 7,19\pmod {20},\\1&\t{otherwise}.\end{cases}
\end{equation}
Also,
\begin{equation}\label{7.5}\prod_{1\ls i<j\ls(p-1)/2\atop p\nmid i^2+3ij+j^2}\l(\f{i^2+3ij+j^2}p\r)
=\begin{cases}-1&\t{if}\ p\eq 19,23,27,31\pmod {40},\\1&\t{otherwise}.\end{cases}
\end{equation}
\end{conjecture}

 Recall that for any prime $p\eq3\pmod4$ the class number $h(-p)$ of the imaginary quadratic field $\Q(\sqrt{-p})$ is odd by \cite{M}.

\begin{conjecture}\label{Conj7.4} Let $\da\in\{\pm1\}$.

{\rm (i)} For any prime $p\eq1\pmod{12}$, we have
\begin{equation}\label{7.6}T_p(1,4\da,1)\eq-3^{(p-1)/4}\pmod p.\end{equation}

{\rm (ii)} Let $p>3$ be a prime. Then
\begin{equation}\label{7.7}\begin{aligned}&\prod_{1\ls i<j\ls(p-1)/2\atop p\nmid i^2+\da4ij+j^2}\l(\f{i^2+\da4ij+j^2}p\r)
\\=&\begin{cases}1&\t{if}\ p\eq1\pmod{24},\\(-1)^{|\{0<k<\f p4:\ (\f kp)=-1\}|}&\t{if}\ p\eq17\pmod{24},
\\\da(-1)^{|\{0<k<\f{p}{12}:\ (\f kp)=-1\}|-1}&\t{if}\ p\eq7\pmod{24},
\\\da(-1)^{|\{0<k<\f{p}{12}:\ (\f kp)=-1\}|+\f{h(-p)-1}2}&\t{if}\ p\eq19\pmod{24}.
\end{cases}\end{aligned}\end{equation}
\end{conjecture}

\begin{conjecture}\label{Conj7.5} Let $p>3$ be a prime. Then
$$\prod_{1\ls i\ls j\ls(p-1)/2\atop p\nmid 4i^2+j^2}\l(\f{4i^2+j^2}p\r)=
\begin{cases} 1&\t{if}\ p\eq1,7,9,19\pmod{20},\\-1&\t{otherwise}.\end{cases}$$
\end{conjecture}

\begin{conjecture}\label{Conj7.6} Let $p>3$ be a prime. Then
$$(-1)^{|\{1\ls k<p/3:\ (\f kp)=-1\}|}\prod_{i,j=1\atop p\nmid 3i+j}^{(p-1)/2}\l(\f{3i+j}p\r)
=\begin{cases}1&\t{if}\ p\eq\pm1\pmod{12},\\(-1)^{\lfloor p/12\rfloor}&\t{if}\ p\eq\pm5\pmod{12},
\end{cases}$$
and
$$\prod_{i,j=1\atop p\nmid 3i-j}^{(p-1)/2}\l(\f{3i-j}p\r)
=\begin{cases}(-1)^{|\{1\ls k<p/3:\ (\f kp)=-1\}|+(p-1)/12}&\t{if}\ p\eq1\pmod{12},
\\(-1)^{|\{1\ls k<p/3:\ (\f kp)=-1\}|-1}&\t{if}\ p\eq5\pmod{12},
\\(-1)^{|\{1\ls k<p/6:\ (\f kp)=1\}|+(p+1)/4}&\t{if}\ p\eq7\pmod{12},
\\-1&\t{if}\ p\eq11\pmod{12}.
\end{cases}.$$
\end{conjecture}

\begin{conjecture}\label{Conj7.7} Let $p>3$ be a prime. Then
$$\prod_{i,j=1\atop p\nmid 4i+j}^{(p-1)/2}\l(\f{4i+j}p\r)
=\begin{cases}1&\t{if}\ p\eq1\pmod{4},\\(-1)^{(h(-p)-1)/2+\lfloor p/8\rfloor}&\t{if}\ p\eq3\pmod{8},
\\(-1)^{|\{1\ls k<p/4:\ (\f kp)=-1\}|}&\t{if}\ p\eq7\pmod{8},
\end{cases}$$
and
$$\prod_{i,j=1\atop p\nmid 4i-j}^{(p-1)/2}\l(\f{4i-j}p\r)
=\begin{cases}(-1)^{(p-1)/4}&\t{if}\ p\eq1\pmod{4},
\\(-1)^{\lfloor p/8\rfloor}&\t{if}\ p\eq3\pmod{4}.
\end{cases}$$
\end{conjecture}

\begin{conjecture}\label{Conj7.8} Let $p>5$ be a prime. Then
\begin{align*}&(-1)^{|\{1\ls k<p/10:\ (\f kp)=-1\}|}\prod_{i,j=1\atop p\nmid 5i+j}^{(p-1)/2}\l(\f{5i+j}p\r)
\\=&\begin{cases}(-1)^{\lfloor (p+1)/10\rfloor}&\t{if}\ p\eq\pm1,\pm3\pmod{20},
\\(-1)^{\lfloor p/20\rfloor}&\t{if}\ p\eq\pm7\pmod{20},
\\(-1)^{\lfloor(p+9)/20\rfloor}&\t{if}\ p\eq\pm9\pmod{20},
\end{cases}
\end{align*}
and
$$\prod_{i,j=1\atop p\nmid 5i-j}^{(p-1)/2}\l(\f{5i-j}p\r)
=\begin{cases}(-1)^{(h(-p)+1)/2}&\t{if}\ p\eq3,7\pmod{20},
\\(-1)^{(p+9)/20}&\t{if}\ p\eq-9\pmod{20},
\\(-1)^{|\{1\ls k<p/10:\ (\f kp)=-1\}|+(h(-p)-1)/2}&\t{if}\ p\eq-1\pmod{20},
\\(-1)^{|\{1\ls k<p/10:\ (\f kp)=-1\}|+\lfloor(p+3)/20\rfloor}&\t{if}\ p\eq1,-3\pmod{20},
\\(-1)^{|\{1\ls k<p/10:\ (\f kp)=-1\}|+\lfloor(p-3)/10\rfloor}&\t{if}\ p\eq -7,9\pmod{20}.
\end{cases}$$
\end{conjecture}

\begin{conjecture}\label{Conj7.9} For any prime $p>3$, we have
$$\prod_{i,j=1\atop p\nmid 6i+j}^{(p-1)/2}\l(\f{6i+j}p\r)
=\begin{cases}(-1)^{|\{1\ls k<p/12:\ (\f kp)=-1\}|}&\t{if}\ p\eq1\pmod{24},
\\(-1)^{|\{\f{p+3}4\ls k\ls\lfloor\f{p+1}3\rfloor:\ (\f kp)=-1\}|}&\t{if}\ p\eq5,-7,-11\pmod{24},
\\(-1)^{(h(-p)+1)/2+\lfloor (p+1)/24\rfloor}&\t{if}\ p\eq -1,-5\pmod{24},
\\(-1)^{\lfloor p/24\rfloor-1}&\t{if}\ p\eq7,11\pmod{24},
\end{cases}$$
and
$$\prod_{i,j=1\atop p\nmid 6i-j}^{(p-1)/2}\l(\f{6i-j}p\r)=(-1)^{|\{\f{p+2}4<k<\f p3:\ (\f kp)=1\}|}.$$
\end{conjecture}

\begin{conjecture}\label{Conj7.10} Let $p>3$ be a prime. Then
$$(-1)^{|\{1\ls k<p/4:\ (\f kp)=-1\}|}\prod_{i,j=1\atop p\nmid 8i+j}^{(p-1)/2}\l(\f{8i+j}p\r)
=\begin{cases}(-1)^{(p+1)/8}&\t{if}\ p\eq-1\pmod{8},
\\1&\t{otherwise},
\end{cases}$$
and
$$\prod_{i,j=1\atop p\nmid 8i-j}^{(p-1)/2}\l(\f{8i-j}p\r)
=\begin{cases}(-1)^{|\{1\ls k<p/4:\ (\f kp)=1\}|}&\t{if}\ p\eq1\pmod{4},
\\(-1)^{(h(-p)+1)/2+(p-3)/8}&\t{if}\ p\eq3\pmod8,
\\-1&\t{if}\ p\eq 7\pmod8.
\end{cases}$$
\end{conjecture}

\Ack. The research is supported by the Natural Science Foundation of China (Grant No. 11971222). The author would like to thank the the referee for helpful comments.


\begin{thebibliography}{99}


\bibitem {BEW} B. C. Berndt, R. J. Evans and K. S. Williams,
 \emph{Gauss and Jacobi Sums}, John Wiley \& Sons, 1998.


\bibitem{B} K. Burde, Eine Verteilungseigenschaft der Legendresymbole,
\emph{J. Number Theory} \textbf{12} (1980), 273--277.

\bibitem{CDE} S. Chowla, B. Dwork and R. J. Evans, On the mod $p^2$ determination of $\bi{(p-1)/2}{(p-1)/4}$,
\emph{J. Number Theory} \textbf{24} (1986), 188--196.

\bibitem{D} L. E. Dickson, \emph{History of the Theory of Numbers},
Vol. III, AMS Chelsea Publ., 1999.

\bibitem{HW} R. H. Hudson and K. S. Williams, Class number formulae of Dirichlet type,
\emph{Math. Comp.} \textbf{39} (1982), 725--732.

\bibitem{IR} K. Ireland and M. Rosen, \emph{A Classical
Introduction to Modern Number Theory}, 2nd Edition, Graduate Texts in
Math., 84, Springer, New York, 1990.

\bibitem{J} M. Jenkins,  Proof of an Arithmetical Theorem leading, by means of Gauss¡¯ fourth demonstration of Legendre¡¯s law of reciprocity, to the extension of that law,
    \emph{Proc. London Math. Soc.} \textbf{2} (1867) 29--32.

\bibitem{M} L. J. Mordell,  The congruence $((p-1)/2)!\eq\pm1\ (\mo\ p)$,
\emph{Amer. Math. Monthly} \textbf{68} (1961) 145--146.

\bibitem{R} H. Rademacher, \emph{Lectures on Elementary Number Theory},
Blaisdell Publishing Company, New York, 1964.

\bibitem{Ri} P. Ribenboim, \emph{The Book of Prime Number Records},
Springer, New York, 1980.

\bibitem{Su} Z.-H. Sun, Values of Lucas sequences modulo primes,
\emph{Rocky Mountain J. Math.} \textbf{33} (2013) 1123--1145.

\bibitem{SS} Z.-H. Sun and Z.-W. Sun, Fibonacci numbers and Fermat's last theorem,
\emph{Acta Arith.} \textbf{60} (1992) 371--388.

\bibitem{S10} Z.-W. Sun, Binomial coefficients, Catalan numbers and Lucas quotients,
\emph{Sci. China Math.} \textbf{53} (2010) 2473--3488.

\bibitem{S19a} Z.-W. Sun, On some determinants with Legendre symbol entries,
\emph{Finite Fields Appl.} \textbf{56} (2019) 285--307.

\bibitem{S19b} Z.-W. Sun, Quadratic residues and related permutations and identities,
\emph{Finite Fields Appl.} \textbf{59} (2019) 246--283.

\bibitem{W} K. S. Williams, On the quadratic residues (mod $p$) in the interval $(0,p/4)$,
\emph{Canad. Math. Bull.} \textbf{26} (1983) 123--124.

\bibitem{WC} K. S. Williams and J. D. Currie, Class numbers and biquadratic reciprocity,
\emph{Canad. J. Math.} \textbf{34} (1982) 969--988.

\end{thebibliography}
 \end{document}